\documentclass[11pt]{article}

\usepackage{epsfig}
\usepackage{graphics}
\usepackage[square, numbers, comma, sort&compress]{natbib}  
\usepackage{verbatim}   
\usepackage{amsmath,amssymb, amsthm}

\usepackage{hyperref}

\usepackage{color}
\usepackage{multirow}
\newtheorem{lem}{\indent{\sc Lemma}}
\newtheorem{prop}{\indent{\sc Proposition}}
\newtheorem{rmk}{\indent{\sc Remark}} 
\newtheorem{thm}{\indent{\sc Theorem}}

\newcommand{\benu}{\begin{enumerate}}
\newcommand{\eenu}{\end{enumerate}}

\newcommand{\beit}{\begin{itemize}}
\newcommand{\eeit}{\end{itemize}}

\newcommand{\be}{\begin{eqnarray}}
\newcommand{\bec}{\begin{center}}
\newcommand{\eec}{\end{center}}
\newcommand{\beo}{\begin{eqnarray*}}
\newcommand{\Bl}{\Bigl(}
\newcommand{\Br}{\Bigr)}

\newcommand{\ee}{\end{eqnarray}}
\newcommand{\eeo}{\end{eqnarray*}}

\newcommand{\f}{\frac}

\newcommand{\fin}{\frac{i}{n}}
\newcommand{\fini}{\frac{i+1}{n}}

\newcommand{\fract}{\frac{1}{t}}

\newcommand{\fpi}{\frac{1}{\sqrt{2 \pi}}}

\newcommand{\ints}{\int \limits_{0}^{s}}
\newcommand{\intt}{\int \limits_{0}^{t}}

\newcommand{\limtin}{\lim \limits_{t \rightarrow \infty}}

\newcommand{\mR}{\mathbb{R}}

\newcommand{\my}{\mathbf{y}}

\newcommand{\n}{\nonumber}

\newcommand{\ov}{\overline}

\newcommand{\R}{\Rightarrow}

\newcommand{\supt}{\sup \limits_{t>0}}
\newcommand{\suptl}{\sup \limits_{t>1}}

\newcommand{\tg}{\textcolor{black}}
\newcommand{\tr}{\textcolor{black}}

\begin{document}

\title{{Some results on diffusion approximation applied to Adaptive MCMC}}
\author{ G K Basak \footnote{Stat-Math Unit, Indian Statistical Institute, 203 B T Road, Kolkata 108. India.}   
 \& Arunangshu Biswas \footnote{Dept. of Statistics, Presidency University, 87/1 College Street, Kolkata 73.India. }}
\date{}
\maketitle

\begin{abstract}
Adaptive Markov Chain Monte Carlo (AMCMC) is a class of MCMC algorithms where the parameters controlling the convergence of the Markov chains are automatically tuned depending on some or all of the previous history of the chain. 
In this situation the transition kernel of the MCMC changes at each iteration and hence proving convergence is not straight forward. In Basak and Biswas \citep{Basak} the authors, 
applying the diffusion approximation procedure to a specially constructed AMCMC with target distribution $\psi(\cdot)$, arrive at a two-dimensional diffusion processes. This continuous time process is relatively easier than its discrete time counterpart.
Although the diffusion in this case is a degenerate one, we show that it satisfies H\"{o}rmander's hypoellipticity condition and consequently has a positive density on its support. 
Using the method of moments we identify the limiting distribution of the $X$-marginal of the diffusion to be the standard Normal density. 
\end{abstract}

\textbf{Keywords and phrases}: Adaptive MCMC, Diffusion approximation, H\"{o}rmander's Hypoelliptic 
conditions, It\^{o}'s Lemma,  MCMC.\\
\textbf{AMS Subject classification: 60J22, 65C05, 65C30, 65C40}

\section {Introduction}
\label{introduction}

Markov Chain Monte Carlo (MCMC) methods are a class of algorithm used to simulate a sample from an arbitrary 
distribution known only upto a constant. One of the algorithms belonging to this class is the Random Walk 
Metropolis-Hastings (RW MH) sampler. The method involves choosing a Markov chain such that the (unique) 
invariant distribution is the target density of interest. This is done by choosing a proposal density, from
which simulating a sample is possible, and then accepting the generated sample with a certain probability 
(called the MH acceptance probability). For more information see \citep{Rosenthal1}. \\

One disadvantage of this method is that the speed of convergence depends on the proposal density. Bad choices of the parameters of the proposal 
(also called the \textit{tuning parameters}) result in very slow convergence of the chain. Consequently, it is of much importance to know what should be 
the optimal choice of the parameters under some criteria. Seminal results in this direction were given for multivariate symmetric Metropolis-Hastings algorithms with a normal proposal in Gelman \textit{et al.} \citep{Gelman}. The target distribution in their case was the product of the marginal densities. Their prescribed value was an acceptance rate of 0.238 where the dimension $d$ of the Markov chain was very large. However, it was shown that this value works good for moderate $d$ as well. 

In another development by Harrio \textit{et al.} \citep{Haario}, the authors proposed the Adaptive
MCMC (AMCMC), where the tuning parameter(s) will be adapted `on-the-fly'. As an example, these values may not be fixed, but a function of the previous sample values. Hence, the proposal density changes at each iteration. This should be done in such a way that the scaling 
constants involved in the proposal density are the best possible choices in some sense. Naturally, the chain loses it Markovian nature and convergence to an invariant distribution can no longer be guaranteed.  

It should be noted that the AMCMC propsed in the literature was discrete time and hence proving convergence required showing that the dicrete time chain converges to stationarity. Such was the approach in Roberts and Rosenthal \citep{Rosenthal2} where the authors gave some sufficient condition for convergence of the chain. These conditions are not necessary and in some cases difficult to verify. 

This paper approaches the problem from a different standpoint. By applying the diffusion approximation scheme we convert the state space variable,together with the tuning parameter (variable) into a continuous time process. Our gain by such an enterprise is that we can then invoke results in the literature for diffusion processes to infer about its invariant distribution whose marginal can then possibly be identified with the 
target distribution of the MCMC. Sometimes this can be done easily when compared to the discrete time setting.

It should be mentioned at this point that the nature of diffusion approximation by Gelman \textit{et al.} \citep{Gelman} is different from our procedure. In our case the limiting diffusion is obtained by taking the limit of the process as the time difference of  successive jumps of the chain converge to zero. In the paper by Gelman \textit{et al.} the diffusion approximation was done by tending the dimension of the state space diverge to infinity and looking at the first co-ordinates of the random vector which is a Markov chain.  

The paper is arranged as follows. Section \ref{definition} contains the definition of the AMCMC and briefly
mentions the diffusion approximation procedure done in \citep{Basak}. 
Section \ref{Main} contains the main result (Theorem \ref{main Theorem}) of this 
paper, i.e., existence of the invariant distribution of the process along with the identification of the target distribution. 
The various subsections of Section \ref{Main} contributes to the proof of Theorem \ref{main Theorem}.
In Section \ref{section tight} we show that the process 
is tight. This combined with the 
hypoelliptic condition in Section \ref{hypoelliptic} shows that the process admits a smooth invariant distribution. 
After establishing moment conditions of the variables under consideration
in Section \ref{momentx} and Section \ref{timeave},
identification of the target distribution is proved in Section \ref{steins}. 
We end with some pointers towards the future direction in Section \ref{conclusion}.

\section{Definitions}
\label{definition}
\setcounter{equation}{0}
We define the AMCMC in such a way that the scaling parameter in the Normal proposal density is a function of whether
the previous sample was accepted or not (ideally it should not depend only on the previous sample but on the 
whole sequence of sample that has been generated, but computations become more extensive in that case). Here 
we formally define our algorithm:   
\begin{enumerate}

\item Select arbitrary $\{X_0,\theta_0 \}\in \mathcal{R} \times [0,\infty)$
 where $\mathcal{R}$ is the state space. Set $n=1$.

\item
Propose a new move say Y where\\
$Y \sim N(X_{n-1}, \theta_{n-1})$.

\item
Accept the new point with probability $\alpha(X_{n-1}, Y)= \min\{1, \frac{\psi(Y)}{\psi(X_{n-1})}\}$.\\
If the point is accepted set $X_{n}=Y ,\ \xi_{i}=1$; else $X_{i}=X_{n-1},\ \xi_{i}=0$.

\item
Set $\theta_{n}=\theta_{n-1} e^{\frac{1}{\sqrt{n}}(\xi_{n} - p)} \ \ p>0$.

\item
Replace $n$ by $n+1$ and go to Step 2.

\end{enumerate}

To apply the diffusion approximation to the AMCMC we define the continuous time process $X_n(t)$
 for all $n \ge 1$ and for all $t>0$ for any target distribution $\psi(\cdot)$:
\be
X_n(0)&=&x_0 \in \mathbf{R}; \n \\
X_n\left(\frac{i+1}{n}\right)&=& X_n(\frac{i}{n}) + \frac{1}{\sqrt{n}} \theta_{n}(\frac{i}{n}) \xi_n(\frac{i+1}{n}) \epsilon_n(\frac{i+1}{n}), \ \ \mbox{i=0, 1, \ldots},\n \\
X_n(t)&=& X_n({\frac{i}{n}}) ,  \ \ \mbox{if $\frac{i}{n} \le t < \frac{i+1}{n}$} \ \ \mbox{for some integer $i$.} \label{X}
\ee
Here, $\xi_n(\frac{i+1}{n})$ conditionally follows the Bernoulli distribution given by:
\beo
P\left(\xi_n(\frac{i+1}{n})=1 \ \Big| \ X_n(\frac{i}{n}), \theta_n(\frac{i}{n}), \epsilon_n(\fini) \right)&=& \min \{\frac{\psi(X_n(\fin)+ \frac{1}{\sqrt{n}}\theta_n(\frac{i}{n})\epsilon_n(\fini))}{\psi(X_n(\frac{i}{n}))},1 \},
\eeo
and $\{\epsilon_n(\frac{i}{n})\}$ are all independent $N(0,1)$ random variables.
The process $\theta_n(t)$ is defined as :
\be
\theta_n(0) &=& \theta_0 \in \mathbf{R^+} \n \\
\theta_n\left({\frac{i+1}{n}}\right) &=& \theta_n\left({\frac{i}{n}}\right) e^{\frac{1}{\sqrt{n}} (\xi_n(\frac{i+1}{n})-p_n({\frac{i}{n}}))}, \ \ \mbox{i=0, 1, \ldots},\n  \\
\mbox{and \ } \theta_n(t)&=& \theta_n(\frac{i}{n}) , \ \ \mbox{if $\frac{i}{n} \le t < \frac{i+1}{n}$ 
for some integer $i$}. \label{theta}
\ee
Where $p_n(\fin) \approx  1-\frac{p}{\sqrt{n}}$ for some $p>0$.

It has been proved in an earlier paper (see \citep{Basak}) that the limiting SDE governing the dynamics of 
the process is the following:

\begin{thm}(from \citep{Basak})
\label{diffusive limits}
 The limit of the process \ \ $\mathbf{Y}_n(t) := \Big(X_n(t), \  \theta_n(t) \Big)'$, where $X_n(t)$ and
 $\theta_n(t)$ is given by (\ref{X}) and (\ref{theta}) respectively, is governed by the SDE:
\be
d \mathbf{Y}_t &=& b(\mathbf{Y}_t)dt+ \sigma(\mathbf{Y}_t)d\mathbf{W}_t, \ \  \mbox{ with } \  \mathbf{Y}_t = (X_t, \theta_t)',
\label{coupled system}
\ee
where,
\beo
b(\mathbf{Y}_t)=\left(  \frac{\theta_t^2}{2} \frac{\psi'(X_t)}{\psi(X_t)} , \ \  \theta_t \left(p - \frac{ \theta_t}{\sqrt{2 \pi}} \frac{|\psi'(X_t)|}{\psi(X_t)} \right) \right)',\\
\eeo
\beo
\sigma(\mathbf{Y_t})=
\left(
\begin{array}{cc}
\theta_t & 0 \\
0 & 0
\end{array}
\right)
\eeo
and $\mathbf{W}_t$ is a two dimensional Wiener process.
\end{thm}

\section {Main result}
\label{Main}
\setcounter{equation}{0}

In this section we concentrate on the case where the target density is standard Normal (i.e., $\psi(x)= \frac{1}{\sqrt{2\pi}} e^{-\frac{x^2}{2}}$). 
Then the SDE takes the form: 
\be
\label{SDE Normal}
d \mathbf{Y}_t &=& b(\mathbf{Y}_t)dt+ \sigma(\mathbf{Y}_t)d\mathbf{W}_t , \ \ \ \mbox{ where, }  \n\\
b(\mathbf{Y}_t)&=&\left(  -\frac{\theta_t^2}{2} X_t , \ \  \theta_t \left(q - \frac{ \theta_t}{\sqrt{2 \pi}} |X_t| \right) \right)'. 
\ee
and $\sigma(\mathbf{Y}_t)$ remains the same. Throughout the section  we assume $\mathbf{Y}_0 = (X_0, \theta_0)'$
is independent of $\{W_t : t \ge 0\}$.

\begin{rmk}
Equation (\ref{SDE Normal}) when written in a more explicit form becomes :
\beo
d X_t &=& -X_t \frac{\theta_t^2}{2} + \theta_t d W_t \\
d \theta_t &=& \theta_t \Bl q - \frac{\theta_t}{\sqrt{2 \pi}} |X_t|\Br dt
\eeo
It resembles that of a coupled Ornstein Uhlenbeck (OU) process with the diffusion coefficient itself  following a logistic equation. One knows that for a standard OU process the N(0,1) distribution is the invariant distribution. In the above case, it is slightly complicated since the diffusion coefficient is not constant. We show that even then the limiting distribution of the diffusion process is Normal. 
\end{rmk}
\begin{rmk}
It will be shown in Lemma \ref{X_t finite} that $E(X_t^2) < \infty, \ \forall t >0$. This implies that $X_t < \infty$ a.s $\forall t$. From the SDE of $\theta_t$ it is shown (see Equation (\ref{theta_t_bound})) that $\theta_t \le \theta_0 e^{qt}$. Consequently $\theta_t < \infty$ almost surely. Therefore the solutions of Equation (\ref{SDE Normal}) is non-explosive. 
\end{rmk}

Here is the main Theorem of this section:
\begin{thm}
\label{main Theorem}
The $X$-marginal of the invariant distribution of (\ref{SDE Normal}) is $N(0,1)$.
\end{thm}

\textbf{Proof:} The proof of the above Theorem is spread over various subsections. In Section 
\ref{section tight} 
we show that the process $(X_t, \eta_t)$ where $\eta_t = 1/ \theta_t$ is tight. This combined with the 
hypoelliptic condition in Section \ref{hypoelliptic} shows that the process admits a invariant distribution. 
The marginal of the invariant distribution is identified as the target distribution in Section \ref{steins}.

\subsection{Tightness of $(X_t, \eta_t)'$}
\label{section tight}
We first state and prove a lemma.
\begin{lem}
\label{martingale lemma}
Fix $T > 0$ and an integer $k \ge 1$. Assume $E(\theta_0^{2k}) < \infty$.
$\int_{0}^{t} \theta_s^k dW_s$ is a martingale with respect to
 $\Bigl\{\mathcal{F}_t=\sigma(X_s, \theta_s ; 0 \le s \le t ), \ 0 \le t \le T \Bigr\}$   and hence for any
$0 \le t \le T$
\beo
E(\int_{0}^{t} \theta_s^k dW_s)&=& 0 .
\eeo

\end{lem}

\textbf{Proof:} 
{It is sufficient to show that 
the local martingale 
$Z_t:=\int_{0}^{t} \theta_s^k dW_s$ is $L_2$-bounded for all $t \le T$. 
So using It\^{o}'s isometry it suffices to show that 
$$ E \Bl \int_{0}^{T} \theta_s^{2k} ds\Br< \infty .$$}
Now,
\be 
d \theta_t &\le& q \theta_t dt \R \theta_t \le \theta_0 e^{qt} \n \\
\R \theta_t^{2k} &\le&\theta_0^{2k} e^{2kqt} \R E \intt \theta_s^{2k}ds \le E(\theta_0^{2k}) \frac{e^{2kqt} -1 }{2kq} < \infty, \label{theta_t_bound}
\ee
for every $t \in [0,T], \ T < \infty.$
$\hfill{\blacksquare}$\\
\subsubsection{Uniform boundedness of moments of $X_t$}
\label{momentx}
We first prove a lemma that will be required in this subsection and elsewhere. Define $F_t = e^{\intt \theta_u^2 du }$ and for any $k \in \mathbb{N},\ C_k := k(1-(2k-1)a)$, where $a>0$ is a constant such that $C_k >0$. 
\begin{lem}
\label{martingale lemma 2}
If $\{X_t\}$ and $\{\theta_t\}$ are solutions to (\ref{SDE Normal}). Fix any $k \in \mathbb{N}$ then
\be
E \Bl F_t^{-C_k} \intt F_u^{C_k} X_u^{2m-1} \theta_u dW_u \Br &=& 0, \ \mbox{\ for any $m \in \{1,2,\ldots,k\}$}, \label{martingale equation}
\ee
where $X_0$ and $\theta_0$ is such that all its moments are finite. 
\end{lem}

\textbf{Proof:} Fix $m \in \{1,2,\ldots,k\}$. Define $\ov{F}_{t,k} := F_t^{-C_k}$. The LHS in (\ref{martingale equation}) is the expectation of $Z_{t,k} (= Z_{t,k}^{(m)}):= \ov{F}_{t,k} Y_{t,k}$ where $Y_{t,k} (= Y_{t,k}^{(m)}) := \intt F_u^{C_k} X_u^{2m-1} \theta_u dW_u$. We show $E(Z_{t,k}) = 0.$ Applying It\^{o}'s lemma to $Z_{t,k}$ we have 
\be
d Z_{t,k} &=& Y_{t,k} d \ov{F}_{t,k} + \ov{F}_{t,k} dY_{t,k} \n \\
&=& -C_k Y_{t,k} \theta_t^2 \ov{F}_{t,k} dt + \ov{F}_{t,k} X_t^{2m-1} \theta_t  F_t^{C_k} dW_t \n \\
&=& -C_k Z_{t,k} \theta_t ^2 dt + X_t^{2m-1} \theta_t dW_t. \label{Z}
\ee
Now, taking $\tilde{Z}_{t,k} = -Z_{t,k}$, yields  
\be
d\tilde{Z}_{t,k} =  C_k Z_{t,k} \theta_t ^2 dt - X_t^{2m-1} \theta_t dW_t = - C_k \tilde{Z}_{t,k} \theta_t ^2 dt + X_t^{2m-1} \theta_t d\tilde{W}_t \label{tilde Z}
\ee
where $\tilde{W}_t = - W_t \stackrel{d}{=} W_t$. From the definition $\tilde{Z}_{0,k} = - Z_{0,k} = 0 = Z_{0,k}$. Comparing the SDE for $Z_{t,k}$ and $\tilde{Z}_{t,k}$ in (\ref{Z}) and (\ref{tilde Z}) we see that they have the same distribution. 
 Therefore $Z_{t,k}$ and $-{Z}_{t,k}$ have the same distribution, which implies that the distribution of $Z_{t,k}$ is symmetric about 0. Now to conclude $E(Z_{t,k})=0$, \ $\forall t \ge 0$ \ we show $Z_{t,k}$ has finite expectation \ $\forall t \ge 0$.
It is sufficient to show that $E(Z_{t,k}^2)< \infty$, \ $\forall t \ge 0$. Now,
\tg{
\beo
Z_{t,k}^2 & = & F_t^{-2C_k} \Bl \intt F_s^{C_k} X_s^{2m-1} \theta_s dW_s  \Br ^2 \le \Bl \intt F_s^{C_k}  X_s^{2m-1} \theta_s dW_s  \Br^2 \n \\
\eeo
a.s, since $F_t^{-2C_k} \le 1$. Therefore, 
\be
E \Bl Z_{t,k}^2\Br &\le& E \Bl \intt F_s^{C_k} X_s^{2m-1} \theta_s  dW_s \Br^2 = E \Bl \intt \underbrace{F_s^{2 C_k} \theta_s}\underbrace{ X_s^{4m-2} \theta_s} ds \Br \n \\
&\le & E \Bl \Bl \intt F_s^{4 C_k} \theta_s^2  ds \Br ^\frac{1}{2} \Bl \intt X_s^{8m-4} \theta_s^2  ds \Br^{\frac{1}{2}}\Br \n \\
&\le&  \sqrt{ E\Bl \intt F_s^{4 C_k} \theta_s^2 ds \Br E \Bl \intt X_s^{8m-4} \theta_s^2 ds \Br  }, \label{Ztksq}
\ee
where the second equality follows from Ito's Isometry and the last two inequalities follow from the Cauchy Schwartz inequality.}
Now for the first expectation in (\ref{Ztksq}) we have
\be 
E \Bl \intt F_s^{4 C_k} \theta_s^2 ds \Br &=& E \Bl \frac{F_t^{4C_k}-1}{4C_k} \Br 
= \frac{1}{4 C_k} E \Bl F_t^{4C_k} - 1 \Br \n \\
&\le& \frac{1}{4C_k} E \Bl e^{4C_k \theta_0^2 (\frac{e^{2qt}-1}{2q})}\Br < \infty , \ \label{first term}
\ee 
{since from (\ref{theta_t_bound}) $\theta_t^2 \le \theta_0^2 e^{2qt}$}.
For the second term in (\ref{Ztksq}) first note that from (\ref{SDE Normal}) and (\ref{theta_t_bound})
\vspace{-10pt}
\be
X_t &=& X_0 - \intt \frac{X_s \theta_s^2}{2} ds + \intt \theta_s dW_s \n \\
\Rightarrow X_t^{8m-4} \theta_t^2 &\le & D_m^2 \theta_0^2 \Bl X_0^{8m-4} + \Bl  \intt \frac{|X_s| \theta_s^2}{2} ds \Br^{8m-4}  
+  \Bl  \intt \theta_s dW_s \Br^{8m-4}  \Br e^{2qt} \n \\
\Rightarrow \intt X_s^{8m-4} \theta_s^2 ds &\le& D_m^2 \theta_0^2 \Bl X_0^{8m-4}\ints e^{2qs} ds + \intt e^{2qs} \Bl \ints \frac{|X_u| \theta_u^2}{2} du \Br^{8m-4}ds \n \\
&+& \intt e^{2qs} \Bl \ints \theta_u dW_u \Br^{8m-4}ds.  \n 
\ee
This implies that 
\be 
\Rightarrow E \Bl \intt X_s^{8m-4} \theta_s^2 ds \Br &\le& D_m \theta_0^2 \Bl E (X_0^{8m-4})\intt e^{2qs} ds  +  E \intt e^{2qs} \Bl \int_{0}^{s} \frac{|X_u| \theta_u^2}{2} du\Br^{8m-4} ds \n \\
&+&  E \intt e^{2qs} \Bl  \int_{0}^{s}\theta_u dW_u \Br^{8m-4} ds \Br \label{eightm}
\ee
for some constant $D_m >0$ that does not depend on  $X_t$.
Clearly the first expectation in the RHS of (\ref{eightm}) is finite $\forall t \ge 0$. \\
For the second expectation in (\ref{eightm}) we proceed as follows. From ( \ref{SDE Normal}) we have the SDE for $\theta_t$ as
\be
d \theta_t &=& \theta_t \Bl q - \frac{|X_t|}{\sqrt{2 \pi}} \theta_t \Br dt = q \theta_t dt - \frac{|X_t|}{\sqrt{2 \pi}} \theta_t^2 dt \n \\
\Rightarrow e^{-qt} d \theta_t - q \theta_t e^{-qt} &=& - e^{-qt} \frac{|X_t|}{\sqrt{2 \pi}} \theta_t^2 dt  \Rightarrow d \Bl \theta_t e^{-qt}\Br = - e^{-qt} \frac{|X_t|}{\sqrt{2 \pi}} \theta_t^2 dt \n\\
\Rightarrow \theta_t e^{-qt} &=& \theta_0 - \intt e^{-qs} \frac{|X_s|}{\sqrt{2\pi}} \theta_s^2 ds \n \\
\Rightarrow \sqrt{\frac{\pi}{2}} \Bl \theta_0 e^{qt} - \theta_t \Br &=& \intt  e^{q(t-s)}\frac{|X_s| \theta_s^2}{2} ds \label{fpi}
\ee

Therefore
\be
\intt \frac{|X_s| \theta_s^2}{2}ds &\le&  \intt e^{q(t-s)} \frac{|X_s| \theta_s^2}{2} ds
\le \sqrt{\frac{\pi}{2}}  \Bl \theta_0 e^{qt} + \theta_t \Br , \n
\ee
from (\ref{fpi}). 

Plugging the value of $\theta_t$ from (\ref{theta_t_bound}) in (\ref{fpi}) we have, 
\be
\intt \frac{|X_s| \theta_s^2}{2}ds &\le& \sqrt{2 \pi} \theta_0 e^{qt}\n \\
\Rightarrow \Bl \intt  \frac{ |X_s| \theta_s^2}{2}ds \Br^{8m-4} &\le& {\Bl \sqrt{2\pi}\Br }^{8m-4}  \Bl \theta_0 e^{qt} \Br^{8m-4} \n \\
\Rightarrow E \intt e^{2qs} \Bl \int_{0}^{s} \frac{|X_u| \theta_u^2}{2} du\Br^{8m-4} ds &\le& (2 \pi)^{4m-2}   (\intt e^{(8m-2)qs}ds)  E(\theta_0^{8m-4}) \n \\
&<& \infty,
\label{second part}  
\ee
for every $t \ge 0$. Hence the second expectation in the RHS of (\ref{eightm}) is also finite $\forall t \ge 0$.

For the third term in the RHS of (\ref{eightm}) let us define $M_s := | \int_{0}^{s} \theta_u dW_u|$ and $M_s^* = \sup \limits_{0<u<s} M_u$. Denoting $[M]_s$ as the quadratic variation process of $M_s$ we have $[M]_s = \ints \theta_u^2 du $. Now, 
\be
E(M_s)^{8m-4} &\le& E(M_s^*)^{8m-4} \le C_mE([M_s]^{4m-2})  \n \\
&=& C_m E \Bl \int_{0}^{s} \theta_u^2 du \Br^{4m-2} \le \Bl  C_m \ints \theta_0^2 e^{2qu} du \Br^{4m-2}, \label{BDG}
\ee
where the second inequality follows from the Burkholder-Davis-Gundy (BDG) inequality and $C_m \in (0, \infty)$ is a constant.
Interchanging the expectation and integrals in the third term of the RHS of (\ref{eightm}) we get 
\be
E \Bl \intt e^{2qs} \Bl \int_{0}^{s} \theta_u d W_u \Br^{8m-4} ds \Br &=&  \intt e^{2qs} E \Bl \int_{0}^{s} \theta_u d W_s \Br^{8m-4} ds  \n \\
&=& \intt e^{2qs} E M_s^{8m-4} ds \n \\
&\le& E(\theta_0^{8m-4} )\intt e^{2qs} \Bl C_m \ints e^{2qu} du \Br^{4m-2}ds \n \\
&<& \infty, \label{third part}
\ee
$\forall t > 0$, where the last but one inequality follows from (\ref{BDG}). Hence the third term of the RHS of (\ref{eightm}) is also finite $\forall t \ge 0$. Hence combining (\ref{second part}) and (\ref{third part}) we have
\beo
E(Z_{t,k}^2) &<& \infty.
\eeo
This combined with the fact that $Z_{t,k}$ is symmetric about zero proves $E(Z_{t,k}) = 0$ and hence the lemma. 
$\hfill{\blacksquare}$ \\
The statement of the above lemma is true even for even powers of $X$, that is

\begin{lem}
\label{rmk mart}
Under the hypothesis of Lemma \ref{martingale lemma 2} the following is true
$$ E \Bl F_t^{-C_k} \intt F_u^{C_k} X_u ^{2m} \theta_u d W_u \Br = 0 \ \mbox{for} \ m \in \{0,1,2,\ldots,k\}.$$
\end{lem}
\textbf{Proof:} We have to prove that $Z_{t,k} = \ov{F_{t,k}} Y_{t,k} := F_t^{-C_k} \intt F_u^{C_k} X_u^{2m} \theta_u dW_u $ has mean zero. Now
\be
d Z_{t,k} &=& -C_k Z_{t,k} \theta_t^2 dt + X_t^{2m} \theta_t d W_t. 
\ee
Define $\ov{Z}_{t,k} = - Z_{t,k}$ and then we see that $Z_{t,k}$ and $-Z_{t,k}$ has the same distribution. We need to show that $Z_{t,k}$ is square integrable. Following steps similar to Equation (\ref{Ztksq}) of the previous lemma 
\be 
E(Z_{t,k}^2) &\le& E \Bl \intt \underbrace{F_s^{2C_k} \theta_s} \underbrace{X_s^{2m} \theta_s}  ds\Br \n \\
&\le& \sqrt{E \Bl \intt F_s^{4C_k} \theta_s^2 ds \Br E \Bl \intt X_s^{8m} \theta_s^2 ds \Br}. \n 
\ee
The first expectation is finite by Equation (\ref{first term}) of Lemma \ref{martingale equation}. For the second expectation we have 
\be
 E \Bl \intt X_s^{8m} \theta_s^2 ds \Br &\le& D_m \theta_0^2 \Bl E(X_0^{8m}) \intt e^{2qs} ds + E \intt e^{2qs} \Bl \intt \frac{|X_u| \theta_u^2}{2} \Br^{8m}\Br ds \n \\
 &+& E \intt e^{2qs} \Bl \intt \theta_u dW_u \Br^{8m} ds 
\ee
By applying methods in the Lemma \ref{martingale equation} the second and the third term can be shown to be finite. This proves the lemma. 
\hfill{$\blacksquare$}\\

Here is the main lemma of this subsection.
\begin{lem}
\label{X_t finite}
For any $k \in \mathbb{N}$,  the $2k^{th}$ ordered moment of $X_t$ is uniformly bounded in $t$, i.e.,
$$ \supt E(X_t^{2k}) < \infty , $$ \tg{if $X_0$ and $\theta_0$ admit finite moments of all order}.
\end{lem}

\textbf{Proof:}
Applying It\^{o}'s lemma to $Y_t = X_t^{2k}$ we get 
\beo
dX_t^{2k} &=& 2k X_t^{2k-1} dX_t + k(2k-1) X_t^{2k-2} \theta_t^2 dt\\
&=& \Bl -kX_t^{2k} \theta_t^2 + k(2k-1) X_t^{2k-2} \theta_t^2 \Br dt + 2k X_t^{2k-1} \theta_t dW_t\\
&\le& \Bl -kX_t^{2k} \theta_t^2 + k(2k-1)(aX_t^{2k} +b) \theta_t^2  \Br dt + 2kX_t^{2k-1} \theta_t dW_t ,
\eeo
since for any fixed $k \in \mathbb{N}$ and small $a > 0$, there exists $b( =b_k)$ large enough such that, $x^{2k-2} < ax^{2k} + b , \ \ \forall x \in \mathbb{R}$.

Thus, for $0 < a < 1/(2k-1)$ we have 
\be
d X_t^{2k} &\le& -X_t^{2k} \theta_t^2 \Bl k - k(2k-1)a \Br dt \n \\
&+& k(2k-1) b \theta_t^2 dt + 2kX_t^{2k-1} \theta_t dW_t\n \\
\Rightarrow d X_t^{2k} + C_k X_t^{2k} \theta_t^2 dt &\le& k(2k-1) b \theta_t^2 dt + 2kX_t^{2k-1} \theta_t dW_t, \label{Xt2k1}
\ee
where $C_k $ and $F_t^{C_k}$ are defined earlier.
Multiplying by the integrating factor $F_t^{C_k}$ on both sides of (\ref{Xt2k1}) we get 
\beo
d \Bl  X_t^{2k} F_t^{C_k} \Br &\le& k(2k-1)b \theta_t^2 F_t^{C_k} dt +  2k  F_t^{C_k} X_t^{2k-1} \theta_t d W_t \\
\Rightarrow X_t^{2k} F_t^{C_k} &\le& X_0^{2k} + k(2k-1)b \intt \theta_u^2 F_u^{C_k}du + 2k \intt F_u^{C_k} X_u^{2k-1} \theta_u dW_u\\
\Rightarrow X_t^{2k} &\le& X_0^{2k}  F_t^{-C_k} + k(2k-1)b F_t^{-C_k} \intt \theta_u^2 F_u^{C_k} du \\
&+& 2k F_t^{-C_k} \intt F_u^{C_k} X_u^{2k-1} \theta_u dW_u .
\eeo
Now,$$ \intt \theta_u^2 F_u^{C_k} du = (F_t^{C_k} -1)/{C_k}$$ 
\be
\Rightarrow E\Bl X_t^{2k} \Br &\le& E \Bl F_t^{-C_k} X_0^{2k} \Br +  k(2k-1) b E \Bl  \frac{1}{C_k} (1 - F_t^{-C_k})\Br \n \\
&+& 2k  E \Bl F_t^{-C_k} \intt F_u^{C_k} X_u^{2k-1} \theta_u dW_u \Br  \label{Xt square}
\ee
For the first term in (\ref{Xt square}) we have that, 
\beo
E(F_t^{-C_k} X_0^{2k}) & \le &  E( X_0^{2k} ) < m< \infty,\ \forall t \ge 0,
\eeo
since $C_k \intt \theta_u^2 du >0.$
Similarly $E \Bl \frac{1}{C_k}(1 - F_t^{-C_k})\Br \le \frac{1}{C_k}$. The third expectation is zero by Lemma \ref{martingale lemma 2}. This proves the lemma.
$\hfill{\blacksquare}$

\subsubsection{Uniform boundedness of moments of $\eta_t=\frac{1}{\theta_t}$}
\begin{lem}
\label{eta}
For any $k \in \mathbb{N}$, the $2k^{th}$ order moment of $\eta_t$ is uniformly bounded in $t \ge 0$, i.e., 
$$ \supt E(\eta_t^{2k}) < \infty , $$ if $X_0$ and $\eta_0$ admit finite moments of all orders.
\end{lem}
\textbf{Proof.}
Take $\eta_t = \frac{1}{\theta_t}$. Then
\beo
d\eta_t &=& - \frac{1}{\theta_t^2}d\theta_t \\
&=& -\frac{1}{\theta_t^2} \theta_t(q - \frac{1}{\sqrt{2\pi}}|X_t| \theta_t)dt 
= - \eta_t (q - \frac{|X_t|}{\eta_t \sqrt{2\pi}})dt 
= (- \eta_t q+ \frac{|X_t|}{\sqrt{2\pi}})dt
\eeo
Multiplying by the integrating factor $e^{qt}$ on both sides of the above equation we get:
\be
d(e^{qt} \eta_t) &=&\frac{e^{qt}|X_t|}{\sqrt{2\pi}}dt \n \\
\Rightarrow e^{qt}\eta_t - \eta_0 &=& \int_{0}^{t}\frac{1}{\sqrt{2\pi}}e^{qu}|X_u|du \n \\
\Rightarrow \eta_t &=& \eta_0 e^{-qt} + \int_{0}^{t} e^{-q(t-u)} \frac{|X_u|}{\sqrt{2\pi}} du  \label{eta form}\\
\Rightarrow E(\eta_t^{2k})  &=& E \Bigl(\eta_0 e^{-qt} + \int_{0}^{t} e^{-q(t-u)} \frac{|X_u|}{\sqrt{2\pi}} du \Bigr)^{2k} \n \\
&\le & 2^{2k-1}\Bigl[ E(\eta_0 e^{-qt})^{2k} + E\Bigl(\int_{0}^{t} e^{-q(t-u)} \frac{|X_u|}{\sqrt{2\pi}} du \Bigr)^{2k} \Bigr] . 
\ee
Now
\beo
\Bigl(\int_{0}^{t} e^{-q(t-u)} \frac{|X_u|}{\sqrt{2\pi}} du \Bigr)^{2k} 
&=& \Bigl(e^{-qt} \int_{0}^{t} e^{qu} \frac{|X_u|}{\sqrt{2\pi}} du \Bigr)^{2k} \\
&=& \frac{(e^{qt} - 1)^{2k}}{(q e^{qt}\sqrt{2 \pi})^{2k}} \Bigl(\frac{q}{e^{qt} - 1} \int_{0}^{t} e^{qu} {|X_u|} du \Bigr)^{2k} \\
&\le & \frac{(e^{qt} - 1)^{2k}}{(q \sqrt{2\pi} e^{qt})^{2k}} \Bigl(\frac{q}{e^{qt} - 1} \int_{0}^{t} e^{qu} |X_u|^{2k} du \Bigr), 
\eeo
where the last inequality follows from the fact that $(E_P(|X|))^{2k} \le E_P(|X|^{2k})$ where $k \in \mathbb{N}$ and $P$ is any probability measure. In the above we take $P(dx) = \frac{q}{e^{qt}-1} e^{qx} dx$ on $[0,t]$. Therefore interchanging the expectation and integrals on the last term of \ref{eta form} we have
\be
\label{eta_bound}
E(\eta_t^{2k}) &\le & 2^{2k-1} \Bigl[ E(\eta_0^{2k}) e^{-2kqt} + \frac{(e^{qt} - 1)^{2k}}{(q \sqrt{2\pi} e^{qt})^{2k}} \frac{q}{e^{qt} - 1} 
\int_{0}^{t} e^{qu} E(|X_u|^{2k}) du \Bigr] \n\\
&\le & 
2^{2k-1} \Bigl[ E(\eta_0^{2k}) + \frac{(e^{qt} - 1)^{2k}}{(q \sqrt{2\pi} e^{qt})^{2k}} M_0 
\Bigr] \n\\
&\le & M_1 < \infty 
\ee
where the last but one inequality follows from Lemma \ref{momentx} that even moments of $X_t$ are uniformly bounded in $t \ge 0$.
$\hfill{\blacksquare}$

\begin{rmk}
 From (\ref{eta_bound}) it is evident that for all $t > 0$, there is a null set, outside of which
$ \theta_t = \frac{1}{n_t} > 0 $ \  whenever  $\theta_0 > 0$, as otherwise, $\supt E(\eta_t^{2k})$ would be infinity.
Again, from the proof above, it is clear that
\beo
\eta_t &=& \eta_0 e^{-qt} + \int_{0}^{t} e^{-q(t-u)} \frac{|X_u|}{\sqrt{2\pi}} du \  >  0, \ \
\mbox{whenever $\eta_0 \ge 0$} .
\eeo
\end{rmk}

Combining the above two lemmas we have the following tightness result for the vector $(X_t, \eta_t)'$.

\subsubsection{Tightness}
\label{tight_xe}

\begin{lem}
\label{lemma tightness}
If $X_0$ and $\theta_0$ admits moments of all orders and $\theta_0 > 0 \ a.s.$then, for the coupled system
(\ref{SDE Normal})
joint distribution of $\{(X_t, \eta_t)' : t \ge 0\}$ is tight.
\end{lem}
\textbf{Proof.}
Let $R_1$ and $R_2$ be two positive numbers. Then
\beo
P(|X_t|<R_1, |\eta_t|< R_2) &=& 1 - P( (|X_t|>R_1)\cup (|\eta_t|>R_2))\\
&>& 1 - (P(|X_t|>R_1) + P(|\eta_t|>R_2))\\
&>& 1- E(|X_t|)/R_1 - E(|\eta_t|)/R_2.
\eeo
Hence given any $\epsilon>0$ we can choose $R_1, R_2$ sufficiently large so that
 $P(|X_t|<R_1, |\eta_t|< R_2) > 1-\epsilon$. This proves the tightness of $(X_t, \eta_t)'$.
 $\hfill{\blacksquare}$

\subsubsection{Finiteness of Time average of moments of $\theta_t$}
\label{timeave}

In this section $C$ will stand for a generic finite constant that might take different values in different 
situations. We assume throughout that $X_0$  and $\theta_0$ admit finite moments of all orders. For non-random initial data this is trivially true.
\begin{lem}
\label{timeave lem}
Let $X_0$ and $\theta_0$ admit finite moments of all order. Then 
\beo
\suptl \fract \intt E(\theta_u^{\frac{k}{2}}) du &<& C \ \mbox{for every} \ k \in \mathbb{N}.
\eeo
\end{lem}
\textbf{Proof: }
We proceed sequentially through the following steps.\\
\underline{Step 1:} We first prove $$\suptl \frac{1}{t} \int_{0}^{t} E(|X_u| \theta_u) du < \infty.$$ 
\vspace{-10pt}
This fact will be used in Step 2. To prove this note that
\be 
d (1 + \theta_t) &=& d \theta_t = \theta_t ( q - |X_t| \theta_t/\sqrt{2 \pi}) dt \n \\
&=& q \theta_t dt - \frac{(1 + \theta_t) |X_t| \theta_t}{\sqrt{2 \pi}} dt + \frac{|X_t| \theta_t}{\sqrt{2 \pi}} dt \n \\
\R d(1 + \theta_t) + \frac{(1 + \theta_t) |X_t| \theta_t}{\sqrt{2 \pi}} dt &=& q \theta_t dt + \frac{|X_t| \theta_t}{\sqrt{2 \pi}} dt \n \\
\R\frac{d(1+\theta_t)}{1+ \theta_t}+ \frac{1}{\sqrt{2\pi}} |X_t| \theta_t dt &=& \frac{\theta_t}{1 + \theta_t} \Bl q + \frac{|X_t|}{\sqrt{2 \pi}}  \Br dt \n \\
&\le& \Bl q + \frac{|X_t|}{\sqrt{2 \pi}} \Br dt \n \\
\Rightarrow \log \frac {1+ \theta_t}{1+ \theta_0} + \fpi \intt |X_u| \theta_u du &\le & qt + \fpi \intt |X_u|du \n \\
\Rightarrow  \frac{1}{t} \intt |X_u| \theta_u du &\le& \sqrt{2 \pi} q+  \frac{1}{t} \intt |X_u|du \n \\
&+& \sqrt{2 \pi} \frac{\log(1 + \theta_0)}{t}. \label{X_theta}
\ee 
Thus, $ \frac{1}{t} \intt E(|X_u| \theta_u) du \le \sqrt{2 \pi} q+  \frac{1}{t} \intt E(|X_u|) du + \sqrt{2 \pi} \frac{E(\log(1 + \theta_0))}{t}$.
Therefore, using the moment bounds for $X_t$ from Section \ref{momentx},
\be 
\sup_{t>1} \frac{1}{t} \intt E(|X_u| \theta_u) du &<& C. \label{E_X_theta}
\ee 
\underline{Step 2:}
We now prove by induction, that for any $k \in \mathbb{N}$,

\be 
\label{theta_moment}
\suptl \frac{1}{t} \intt  E(\theta_u^{\frac{k}{2}}) du 
< C .
\ee

Let, as before, $\eta_t = \frac{1}{\theta_t}$ then $d \eta_t =( - q \eta_t + |X_u|/ \sqrt{2\pi}) dt $. \\

Applying It\^{o}'s lemma to $Y_t = X_t^2 \eta_t^{2 - k/2}$, with $k \in \mathbb{N}$, we get
\be
dY_t &=& 2 X_t \eta_t^{2 - k/2} dX_t + (2-k/2) X_t ^2 \eta_t^{1-k/2} d \eta_t + \frac{1}{2} 2 \eta_t^{2 - k/2} (dX_t)^2 \n \\
&=& 2 X_t \eta_t^{2 - k/2} (- \frac{X_t}{2 \eta_t^2} dt + \frac{1}{\eta_t} dW_t) + (2 - k/2) X_t^2 \eta_t^{1 - k/2} ( - q \eta_t dt + \frac{|X_t|} {\sqrt{2 \pi} }dt ) \n \\
&+&  \eta_t^{2 - k/2} \eta_t^{-2} dt  \n \\
&=& \Bigl( - X_t^2  \eta_t^{- k/2} - q (2 - k/2) X_t^2 \eta_t^{2 - k/2} + \frac{2 - k/2}{\sqrt{2 \pi}} |X_t|^{3} \eta_t^{1 - k/2} + \eta_t^{- k/2} \Bigr) dt \n \\
&+& 2 X_t \eta_t^{1 - k/2} dW_t . \label{stoc eqn}
\ee
Thus, integrating both side from $0$ to $t$, rearranging and dividing by $t$ and then taking expectations we get
\be
\intt \theta_s^{\frac{k}{2}} ds &=&  X_t^2 \eta_t^{\frac{4-k}{2}} - X_0^2 \eta_0^{\frac{4-k}{2}} + \intt X_s^2 \theta_s^{\frac{k}{2}} ds \n \\
&+& \frac{(4-k)q}{2} \intt X_s^2 \eta_s^{\frac{4-k}{2}} ds 
- \frac{2 - k/2}{\sqrt{2 \pi}} \intt |X_s|^3 \eta_s^{\frac{2-k}{2}} ds \n \\
&-& 2 \intt X_s \eta_s^{\frac{2-k}{2}} dW_s \n \\
\Rightarrow \frac{1}{t} \intt \theta_s^{\frac{k}{2}} ds  &=& \frac{1}{t}(X_t^2 \eta_t^{\frac{4-k}{2}} - X_0^2 \eta_0^{\frac{4-k}{2}}) + \frac{1}{t} \intt X_s^2 \theta_s^{\frac{k}{2}} ds \n \\
&+&  \frac{(4-k)q}{2t} \intt X_s^2 \eta_s^{\frac{4-k}{2}} ds \n \\
&& - \frac{2 - k/2}{t\sqrt{2 \pi}} \intt |X_s|^3 \eta_s^{\frac{2-k}{2}} ds - \frac{2}{t} \intt X_s \eta_s^{\frac{2-k}{2}} dW_s \n \\
\Rightarrow \suptl \frac{1}{t} \intt E(\theta_s^{\frac{k}{2}}) ds &\le&  \suptl \frac{1}{t}E (X_t^2 \eta_t^{\frac{4-k}{2}} - X_0^2 \eta_0^{\frac{4-k}{2}}) + \suptl \frac{1}{t} \intt E(X_s^2 \theta_s^{\frac{k}{2}}) ds \n \\
&+& \frac{(4-k)q}{2} \suptl \fract \intt E(X_s^2 \eta_s^{\frac{4-k}{2}}) ds \n \\
&-& \frac{4-k}{2 \sqrt{2 \pi}}  \suptl \fract \intt E(|X_s|^3 \eta_s^{\frac{2-k}{2}}) \n \\
&-&  2 \suptl \fract E \intt (X_s \eta_s^{\frac{2-k}{2}}) dW_s
\label{int_theta1}.
\ee
Now for any $k \in \mathbb{N}$ we have, 
\be
\frac{1}{t} \intt X_s^2 \theta_s^{\frac{k}{2}} ds &=& \frac{1}{t} \intt (|X_s|^{\frac{k}{k+1}} 
\theta_s^{\frac{k}{2}}) ( |X_s|^{\frac{k+2}{k+1}} ) ds  \n \\
&\le&   \Bl \frac{1}{t} \intt  |X_s| \theta_s^{\frac{k+1}{2}} ds \Br ^{ \frac{k}{k+1}} 
\Bl \frac{1}{t} \intt |X_s|^{k+2} ds \Br ^ {\frac{1}{k+1}},  
\ee
which follows from the Holder's inequality with $p = \frac{k+1}{k}$ and $q = k+1$. Therefore, 
\be
E \Bl \fract \intt X_s^2 \theta_s^{\frac{k}{2}} \Br 
&\le& E \Bl \Bl \frac{1}{t} \intt  |X_s| \theta_s^{\frac{k+1}{2}} ds \Br ^{ \frac{k}{k+1}} 
\Bl \frac{1}{t} \intt |X_s|^{k+2} ds \Br ^ {\frac{1}{k+1}} \Br \n\\
&\le &   \Bl E \Bl \frac{1}{t} \intt  |X_s| \theta_s^{\frac{k+1}{2}} ds \Br \Br ^{ \frac{k}{k+1}}\times \Bl E \Bl \frac{1}{t} \intt |X_s|^{k+2} ds \Br \Br ^ {\frac{1}{k+1}}  \n\\
&=& \Bl \frac{1}{t} \intt E(|X_s| \theta_s^{\frac{k+1}{2}}) ds \Br ^{ \frac{k}{k+1}}
 \Bl \frac{1}{t} \intt E(|X_s|^{k+2}) ds  \Br ^ {\frac{1}{k+1}},
 \ee 
 where the last inequality follows from Holder's inequality with $p = \frac{k+1}{k}$ and $q = k+1$. Therefore,
 \be
\suptl \fract \intt E(X_s^2 \theta_s^{\frac{k}{2}}) ds &\le& 
 \Bl \suptl \frac{1}{t} \intt  E(|X_s| \theta_s^{\frac{k+1}{2}}) ds  \Br ^{ \frac{k}{k+1}} \n \\
 &\times & \Bl \suptl \frac{1}{t} \intt E(|X_s|^{k+2}) ds \Br ^ {\frac{1}{k+1}}. \label{prefourth}
\ee
\vspace{-10pt}
Again $\forall k \ge 2$, 
\be
\label{induction}
d\theta_t^{\frac{k-1}{2}} = \frac{k-1}{2} \theta_t^{\frac{k-1}{2} - 1} d\theta_t 
&=& \frac{k-1}{2} \theta_t^{\frac{k-1}{2}} ( q - \frac{|X_t| \theta_t}{\sqrt{2 \pi}}) dt \n\\
&=& \frac{q(k-1)}{2}  \theta_t^{\frac{k-1}{2}} dt  - \frac{k-1}{2} \frac{|X_t| \theta_t^{\frac{k+1}{2}}}{\sqrt{2 \pi}} dt \n\\
\Rightarrow \fpi \intt |X_s| \theta_s^{\frac{k+1}{2}} ds 
&=& q \intt \theta_s^{\frac{k-1}{2}} ds - \frac{2}{k-1} (\theta_t^{\frac{k-1}{2}} - \theta_0^{\frac{k-1}{2}}) \n\\
\Rightarrow \suptl \fract \intt E \Bl |X_s| \theta_s^{\frac{k+1}{2}} ds \Br 
&\le& \sqrt{2 \pi} q \suptl \fract \intt E \Bl \theta_s^{\frac{k-1}{2}} ds \Br +  \frac{2 \sqrt{2 \pi}}{k-1} \suptl \fract  E \Bl \theta_0^{\frac{k-1}{2}}\Br \n \\\label{X_s_theta_kp1}.
\ee
Plugging (\ref{X_s_theta_kp1}) in (\ref{prefourth}) 
\be
\suptl \fract \intt E(X_s^2 \theta_s^{\frac{k}{2}}) ds &\le& \Bigl(\sqrt{2 \pi}q\suptl \fract \intt E \Bl \theta_s^{\frac{k-1}{2}} \Br ds + \frac{2 \sqrt{2 \pi}}{k-1} \suptl \fract E \Bl \theta_0^{\frac{k-1}{2}}\Br \Bigr)^\frac{k}{k+1} \n \\
& \times &  \Bl \suptl \frac{1}{t} \intt E(|X_s|^{k+2}) ds \Br ^ {\frac{1}{k+1}}  \label{X_s_square}.
\ee

And finally plugging (\ref{X_s_square}) in (\ref{int_theta1})  we get for $k \ge 2 $ 
\be
\suptl \frac{1}{t} \intt E(\theta_s^{\frac{k}{2}}) ds &\le&  \suptl \frac{1}{t}E (X_t^2 \eta_t^{\frac{4-k}{2}} - X_0^2 \eta_0^{\frac{4-k}{2}}) \n \\
&+& \Bigl( \sqrt{2 \pi} q \suptl \fract \intt E \Bl \theta_s^{\frac{k-1}{2}} \Br ds + 
\frac{2\sqrt{2 \pi}}{k-1} \suptl \fract E ( \theta_0^{\frac{k-1}{2}}) \Bigr)^\frac{k}{k+1} \n \\
&\times & \Bl \suptl \frac{1}{t} \intt E(|X_s|^{k+2}) ds \Br ^ {\frac{1}{k+1}} \n \\
&+&  \frac{(4-k)q}{2} \suptl \fract \intt E(X_s^2 \eta_s^{\frac{4-k}{2}}) ds \n \\
&-& \frac{4-k}{2 \sqrt{2 \pi}}  \suptl \fract \intt E(|X_s|^3 \eta_s^{\frac{2-k}{2}})ds \n \\
&-&  2 \suptl \fract E \intt (X_s \eta_s^{\frac{2-k}{2}}) dW_s \label{int_theta}.
\ee

To prove $\suptl \fract \intt E(\theta_s^\frac{k}{2})ds$
is finite $\forall \ k \in \mathbb{N}$ we proceed by induction:\\
\underline{Step 2a}: For k=1 we consider Equation (\ref{int_theta1}). By an application of the Young's  inequality and the fact that all the moments of $X_s$ and $\eta_s$ are uniformly bounded (proved earlier in Lemma \ref{X_t finite} and \ref{eta}) we have:
\beo
\suptl \frac{1}{t} E(X_t^2 \eta_t^{\frac{4-1}{2}} - X_0^2 \eta_0^{\frac{4-1}{2}}) &<& C ,
\label{first_k_1}  \\
\suptl \frac{1}{t} \intt E(X_s^2 \eta_s^{\frac{4-1}{2}}) ds 
 &<& C , \label{third_k_1} \\ 
\mbox{and} & & \n \\
\suptl \fract \intt E(|X_s|^3 \eta_s^{\frac{2-1}{2}}) ds &<& C . \label{fourth_k_1}
\eeo
This proves that the first, third and fourth term in the RHS of (\ref{int_theta1}) is finite. The second term of (\ref{int_theta1}) is bounded by the RHS of (\ref{prefourth}), whose first term is finite by (\ref{E_X_theta}) of Step 1 and the second term is finite by the uniform boundedness of moments of $X$. Therefore we are left with only the It\^{o} integral or the last term of (\ref{int_theta1}). Now, $$ E \Bl \intt X_s \eta_s^{\frac{1}{2}}dW_s \Br^2 = E \intt X_s^2 \eta_s ds $$ is finite $\forall t \ge 0$ by an application of Young's inequality and the uniform boundedness of all the moments of $X_t$ and $\eta_t$. Therefore $\intt X_s \eta_s^{\frac{1}{2}} dW_s$ is a square integrable martingale and hence 
$$\suptl \fract E \intt X_s \eta_s^{\frac{1}{2}} dW_s = 0.$$
This completes the proof that $\suptl \fract \intt E(\theta_s^{\frac{1}{2}}) ds$ is finite $\forall t >0$.

\underline{Step 2b:} Assume that the hypothesis is true for $k \le m-1$, for $m \ge 2$ i.e.,
\beo
\suptl \frac{1}{t} \intt E(\theta_s^{\frac{k}{2}}) ds  &<& C,  \ \ {k \le m-1}. \label{ind_hypo}
\eeo  
 
\underline{Step 2c:} Consider $k = m \ge 2.$ In this case we consider Equation (\ref{int_theta}). \\
For $m = 2$ the RHS of (\ref{int_theta}) is finite by the moment bounds of $X_s$ and $\eta_s$ and by the proof that $\suptl \intt E(\theta_s^{\frac{1}{2}})ds < \infty $ in Step 2a.

For $m = 3, 4$ the first term in the RHS of (\ref{int_theta}) is finite (by the arguments given in 2a). The second (product) term is finite by the induction hypothesis (in 2b) and by the finiteness of the moments of $X_s$. The third term is finite by the finiteness of the moments of $X_t$ and $\eta_t$. The fourth term is negative for $m=3$ or zero for $m=4$. Hence it is bounded by zero.

For the fifth (It\^{o} Integral) term in (\ref{int_theta}) we first apply the It\^{o}'s lemma and then Cauchy Schwartz inequality to get 
\be
\intt E(X_s^2 \theta_s^{m-2})  ds \le \intt \sqrt{E(X_s^4)  E(\theta_s^{2(m-2)})}ds &\le& C \intt \sqrt{E(\theta_s^{2(m-2)})} ds \n  \\
&<& \infty \label{fifth term} ,
\ee
since $\theta_t$ is bounded as in Equation (\ref{theta_t_bound}).  
Thus,
$ \intt X_s \theta_s^{\frac{m-2}{2}} d W_s$ is a square integrable martingale with respect to the given filtration over any finite interval $[0,T]$ and therefore the expectation is zero.\\
Next consider $m >4$. 
For the first term in the RHS of (\ref{int_theta}) apply Young's inequality with $p = m-3$ and $q = \frac{m-3}{m-4}$ to get 
\beo
X_s^2 \theta_s^{\frac{m-4}{2}} \le \frac{1}{m-3} X_s^{{2(m-3)}} + \frac{m-4}{m-3} \theta_s^{\frac{m-3}{2}} &\Rightarrow& E(X_s^2 \theta_s^{\frac{m-4}{2}}) \le \frac{1}{m-3} E(X_s^{{2(m-3)}})\\
&+& \frac{m-4}{m-3} E(\theta_s^{\frac{m-3}{2}}) \n \\
\Rightarrow \suptl \fract \intt E (X_s^2 \theta_s^{\frac{m-4}{2}}) ds &\le& \frac{1}{m-3} \suptl \fract \intt E(X_s^{2(m-3)})ds \n \\
&+& \frac{m-4}{m-3} \suptl \fract \intt E(\theta_s^{\frac{m-3}{2}})ds \n \\
&<& \infty,
\eeo
which follows from the fact that moments of $X_t$ are uniformly bounded and the second term is finite by the induction hypothesis. Consequently, the first term in the RHS of (\ref{int_theta}) is finite.\\
The second (product) term is finite by the induction hypothesis and by the finiteness of the moments of $X_s$ (as argued in the case $m = 3,4$ above). \\
The third term is negative.\\ 
The fourth term we apply the Young's inequality  with $p = m-1$ and $q = \frac{m-1}{m-2}$ to get: 
\be
|X_s|^3 \theta_s^{\frac{m-2}{2}} \le \frac{|X_s|^{3p}}{p} + \frac{\theta_s^{\frac{q(m-2)}{2}}}{q} 
&=& \frac{1}{m-1} |X_s|^{3(m-1)} + \frac{m-2}{m-1} \theta_s^{\frac{m-1}{2}} , \n \\
\Rightarrow \suptl \fract \intt E\Bl |X_s|^3 \theta_s^{\frac{m-2}{2}}\Br ds &\le& \frac{1}{m-1} \suptl \fract \intt E \Bl |X_s|^{3(m-1)} \Br ds \n \\
&& + \frac{m-2}{m-1} \suptl \fract \intt  E \Bl |\theta_s|^{\frac{m-1}{2}} \Br ds \n \\
&<& \infty,
\label{X_s_cube}
\ee
which follows from the fact that the moments of $X_t$ are uniformly bounded in $t$ and by the induction hypothesis.\\
For the fifth term we argue as in (\ref{fifth term}) to infer that it is a square integrable martingale with respect to the given filtration over any finite interval $[0,T]$ and hence the expectation is zero. \\
Therefore the LHS of (\ref{int_theta}) is finite for all $m \ge 2$.
Thus the Steps 2a, 2b and 2c complete the proof of Step 2 (\ref{theta_moment}) and therefore Lemma \ref{timeave lem} is proved. $\hfill{\blacksquare}$
 
\subsection{Hypoelliptic condition}
\label{hypoelliptic}

Here we show that the vector fields corresponding to (\ref{SDE Normal}) satisfies the H\"{o}rmander's hypoelliptic 
condition (see the proposition for the statement of the condition).
{Since the condition requires smooth vector fields, we convert the drift and diffusion coefficients in (\ref{coupled system}) into smooth vector fields. \\
For this purpose, define $$b_\epsilon(x, \eta)= \left( -\f {x}{2 \eta^2 } , \   - q\eta + \f{g_{\epsilon}(x)}{\sqrt{2 \pi}} \right),$$ where $g_\epsilon(x),\ \mbox{a smooth function} \ \rightarrow |x|$ as $\epsilon \downarrow 0$ in the  point-wise limits and $\sigma(x, \eta)=\left(\begin{array}{ll}
 1/\eta & 0 \\
0 & 0 \\
\end{array}
\right)$
as the drift and the diffusion coefficient respectively of the equation with the re-parametrisation $\eta = \frac{1}{\theta}$. Such function $g_\epsilon(\cdot)$ can be constructed by convoluting the function $|x|$ with a mollifier (for example $ \frac{1}{\sqrt{2\pi} \epsilon} e^{- \frac{1}{2 \epsilon^2} x^2}$). \\
Consider an SDE in the Stratonovich  form:
\be
dX_t &=& A_0(X_t)dt + \sum \limits_{\alpha=1}^n A_{\alpha}(X_t)\circ dW^{\alpha}_t  \label{SDE_vector}.
\ee
where $A_0, \{A_{\alpha}: \alpha = 1, \ldots, n\}$ is a smooth vector fields on a differential manifold $M$ 
and $\circ$ denotes Stratonovich integral. The SDE in the It\^{o} form and the Stratonovich form are 
interchangeable.
For a multidimensional SDE, given in the It\^{o}'s form,
\beo
d \mathbf{X}_t &=& \mathbf{b}(t,\mathbf{X}_t) d t + \mathbf{\sigma}(t, \mathbf{X}_t) d \mathbf{W}_t
\eeo
can be readily converted into the Stratonovich form from the following equation:
 \beo
 \tilde{b}_i(t,\mathbf{x}) &=& b_i(t, \mathbf{x}) -  \frac{1}{2} \sum \limits_{j=1}^{p} \sum \limits_{k=1}^{n}  \frac{\partial \mathbf{\sigma}_{i,j}}{\partial x_k} \sigma_{k,j}; \; \; 1 \le i \le n
 \eeo
where $\mathbf{\tilde{b}}(t,x)=(\tilde{b}_i(t,x))'$ is the drift term for the Stratonovich form. In our case,
$p=n=2$ and from the form of $\sigma$ in (\ref{SDE Normal}),  we find that $\tilde{b^{\epsilon}}$ and 
$b^{\epsilon}$ are the same and it equals $A_0$. 
We identify the diffusion coefficients $A_1(\mathbf{X}_t) = (\eta , 0)'$ and $A_2(\mathbf{X}_t) = (0 , 0)'$  
as vector fields in $M$, here upper half plane of $\mR^2$.
Here is the condition due to H\"{o}rmander \citep{hormander}:
\begin{prop}
\label{hypo prop}
Let $\{A_0, A_1, \ldots, A_n\}$ be $n+1$ smooth vector fields on  
a smooth manifold $M$. 
Define the Lie Bracket $[V, W]$ 
between two vector fields $V$ and $W$ as another vector field on $M$  defined in the following manner  
$$ [V, W](f)= V(W(f)) - W(V(f)) \ \ \forall f \in C^{\infty}(M).$$
The H\"{o}rmander's hypoelliptic condition is satisfied if :
\beo
A_{j_0}(\my),  \ [A_{j_0}(\my), A_{j_1}(\my)], \ [[A_{j_0}(\my), A_{j_1}(\my)],A_{j_2}(\my)], \\
\ldots  \;
[ [ [ [ A_{j_0}(\my), A_{j_1}(\my) ], A_{j_2}(\my) ], A_{j_3}(\my) ],\ldots, A_{j_k}(\my) ] 
\eeo
spans  
$M$
for every  
$\my \in M$
 and any $1 \le j_0 \le n$ and $\{j_1, \ldots, j_k\} \in \{0, 1, \ldots, n\} , \ k \ge 1$.
 \end{prop}
\begin{lem}
\label{hypo condn}
The vector fields $A_0^{\epsilon}(\my)$ and $A_1(\my)$ satisfy H\"{o}rmander's hypoelliptic condition of Proposition
\ref{hypo prop}.
\end{lem}

\textbf{Proof:}  Identifying (\ref{SDE_vector}) with (\ref{SDE Normal}) we have (writing $\my = (x,\eta)'$):
\beo
A_0^\epsilon(\my) &=& - \frac{x}{2 \eta^2} \frac{\partial}{\partial x} + (-q \eta + \frac{g_{\epsilon}(x)}{\sqrt{2 \pi}}) \frac{\partial}{ \partial \eta}, \\
A_1(\my)&=& \frac{1}{\eta} \frac{\partial}{\partial x}.
\eeo
Therefore the vectors corresponding to  
$A_1(\my)$ and $[A_1(\my), A_0^\epsilon(\my)]$ will be 
$\Bl \frac{1}{\eta}, 0 \Br^{T}$ and $\Bl   \frac{1}{\eta^2} \Bl - \frac{1}{2\eta} - q \eta + \fpi g_\epsilon(x) \Br,  \fpi \frac{1}{\eta} g'_\epsilon(x)\Br^T$. Note, $\theta_t = 1/\eta_t > 0$ almost surely, since by Lemma \ref{eta} we have $\sup \limits_{t>0} E(\eta_t^{2}) < \infty$. Thus, the zero set of $\{X_t\}$ has Lebesgue measure zero almost surely since the zero set of $\{W_t\}$ has Lebesgue measure zero. Therefore these two vector fields span the upper half plane of $\mathbb{R}^2$, for $x \ne 0$.
Also, for $x \ne 0$, we can take $\epsilon \to 0$ and get the same result. Note that the convergence is uniform over each compacts in the set $\{(x, \eta): x \neq 0, \eta >0 \}.$
 $\hfill{\blacksquare}$
 \begin{rmk}
 \label{normal mollifier}
 In the case of the Normal mollifier i.e,
\beo
g_\epsilon(y)&=& \frac{1}{\sqrt{2 \pi} \epsilon}\int_{-\infty}^{\infty}|x|e^{-\frac{1}{2\epsilon^2}(y-x)^2} dx \n \\
&=& \frac{1}{\sqrt{2 \pi} \epsilon}\int_{0}^{\infty}x e^{-\frac{1}{2\epsilon^2}(y-x)^2} dx + \frac{1}{\sqrt{2 \pi} \epsilon}\int_{-\infty}^{0}(-x) e^{-\frac{1}{2\epsilon^2}(y-x)^2} dx. \n  
\eeo 
For the first integral 
\beo
\frac{1}{\sqrt{2 \pi} \epsilon}\int \limits_{0}^{\infty}x e^{-\frac{1}{2\epsilon^2}(y-x)^2} dx &=& \frac{1}{\sqrt{2 \pi}} \int \limits_{-\frac{y}{\epsilon}}^{\infty}(y + \epsilon z) e^{-\frac{z^2}{2}} dz, \ \mbox{substituting $z = \frac{x-y}{\epsilon}$,} \n \\
&=& y (1- \Phi(- \frac{y}{\epsilon})) + \fpi \epsilon \int \limits_{\frac{y^2}{2 \epsilon^2}}^{\infty} e^{-t}dt, \ \mbox{substituting $t = \frac{z^2}{2}$,} \n \\
&=& y \Phi(\frac{y}{\epsilon}) + \fpi \epsilon  e^{-\frac{y^2}{2 \epsilon^2}},
\eeo
where $\Phi(\cdot)$ is the distribution function of the standard Normal variable. Similarly for the second integral we have 
 \beo
 \frac{1}{\sqrt{2 \pi} \epsilon}\int \limits_{-\infty}^{0}(-x) e^{-\frac{1}{2\epsilon^2}(y-x)^2} dx &=& \frac{1}{\sqrt{2 \pi}} \int \limits_{-\infty}^{-\frac{y}{\epsilon}}- \ov{y + \epsilon z} e^{-\frac{z^2}{2}} dz, \ \mbox{substituting $z = \frac{x-y}{\epsilon}$,} \n \\
 &=& - y \Phi(- \frac{y}{\epsilon}) - \fpi \epsilon \int_{-\infty}^{- \frac{y}{\epsilon}}  z e^{-\frac{z^2}{2}} dz \n \\
&=&  - y \Phi(- \frac{y}{\epsilon}) - \epsilon \fpi \int_{\infty}^{\frac{y^2}{2 \epsilon^2}} e^{-t} dt \n \\
&=& - y \Phi(- \frac{y}{\epsilon}) + \epsilon \fpi e^{- \frac{y^2}{2 \epsilon^2}}\n \\
 \R g_\epsilon(y) &=& y \Bl \Phi(\frac{y}{\epsilon}) - \Phi(- \frac{y}{\epsilon}) \Br + 2 \epsilon \phi(\frac{y}{\epsilon}) \n \\
 \R \frac{d}{dy} g_\epsilon'(y) &=& \Phi(\frac{y}{\epsilon}) - \Phi(- \frac{y}{\epsilon}) + 2 \frac{y}{\epsilon} \phi(\frac{y}{\epsilon}) -2 \frac{y}{\epsilon} \phi(\frac{y}{\epsilon})\n \\
  \R |\frac{d}{dy} g_\epsilon(y)|&\le& |\Phi(\frac{y}{\epsilon}) - \Phi(-\frac{y}{\epsilon})|, \n  
  \eeo
{where $\phi(\cdot)$  is the density function of the standard Normal variable.}  Now for any $\epsilon>0$ and any $y \in \mR$ we have 
\beo
|\Phi(\frac{y}{\epsilon}) - \Phi(-\frac{y}{\epsilon})| &\le& 1 \n 
\eeo
which implies that 
\beo
|\frac{d}{dy} g_\epsilon (y)| &\le& 1 \ \forall \epsilon >0 \n \\
\R  \sup_{\epsilon} |\frac{d}{dy} g_\epsilon(y)| &<& \infty,\ \forall \ y \in \mR,
\eeo
which implies that the family $\{ g_\epsilon (\cdot) \}$ is equicontinuous. 
\end{rmk}
\bigskip
It is well known that if the vector fields $A_0(\my)$ and $A_1(\my)$ satisfy the above conditions then the 
solution of the SDE (\ref{SDE_vector}) admits a smooth transition density (see, for example Nualart \citep{nualart}).\\
Hence, even though the original diffusion is singular
its transition probability has density (see Kliemann \citep{kliemann}). Again, since the coupled diffusion is tight,
it admits unique invariant probability by Kliemann \citep{kliemann}
which admits a density.
\begin{rmk}
Note that although we are interested in the distribution of $\{X_t\}$ showing tightness of the process $\{X_t\}$ only it would
 not suffice since $\theta_t$ may be a function of $\{X_s; 0 \le s \le t\}$, so marginally $\{X_t\}$ may not be a 
Markov process. Hence $\sup \limits_{t >0} E|X_t|< M$ would give the tightness of X but it would not be possible to say 
anything about the existence of a unique invariant distribution of $\{X_t\}$.
\end{rmk}

\subsection{Identifying the limiting distribution}
\label{steins}

We first prove a lemma that will be required in this subsection.
For any $s>0$, define $F_s(t) = s \intt \theta_u^2 du. $
\begin{lem}
\label{Fkt}
$$ \limtin E(e^{-F_s(t)}) = 0, \ \forall s >0.$$
\end{lem}

\textbf{Proof:} We prove for $s=1$. The  proof can be carried out in a similar fashion  for any $ s > 0 $.
\be
\frac{F_1(t)}{t} &=&\fract \intt \theta_s^2 ds \ge \frac{1}{\fract \intt \frac{1}{\theta_s^2}ds}  = \frac{1}{\fract \intt \eta_s^2 ds }, \label{F1}
\ee
where the last but one inequality follows from Jensen's (by taking $\psi(x) = \frac{1}{x},\ x > 0 $ which is convex). This implies
\be
\frac{1}{\fract F_1(t)} &\le& \fract \intt \eta_s^2 ds 
\Rightarrow  \frac{1}{F_1(t)} \le \fract \fract \intt \eta_s^2 ds. \label{F_1} 
\ee
Therefore,
\beo
e^{-F_1(t)} &=& \frac{1}{e^{F_1(t)}} \le \frac{1}{{F_1(t)}} \ \mbox{(since $e^x \ge x,\ \forall x >0$)} \n \\
&\le& \fract \fract \intt \eta_s^2 ds \ \ \mbox{(from (\ref{F_1}))}\\
\Rightarrow E(e^{-F_1(t)}) &\le& \fract E\Bl \fract \intt \eta_s^2 ds\Br \le \fract C .
\eeo
where $C = \supt E(\fract \intt \eta_s^2 ds) < \infty$, from Lemma \ref{eta}.
So  
\beo
\limtin E(e^{-F_1(t)}) &=& 0.
\eeo
$\hfill{\blacksquare}$

\begin{lem}
\label{even odd}
Assuming that all the moments of $X_0$ and $\theta_0$ exists we have 
\begin{center}
$
\tr{\lim \limits_{t \to \infty} E(X_t^{r})} = \left \{
\begin{array}{ll}
\frac{(2k)!}{2^k k!} & \mbox{when $r = 2k$}\\
0 & \mbox{when $r = 2k+1$} .
\end{array}
\right.
$
\end{center}

\end{lem}

\textbf{Proof:} We prove using induction for both even and odd moments:

\textbf{Even moments}

\begin{enumerate}
\item We first show $\lim_{t \rightarrow \infty} E(X_t^{2}) = 1$.\\
Applying It\^{o}'s lemma to $X_t^2$ we have 
\beo
dX_t^{2}&=& \Bl -X_t^2 \theta_t^2 + \theta_t^2 \Br dt + 2 X_t \theta_t dW_t .
\eeo
Multiplying by the integrating factor
$e^{F_1(t)}$, where $F_1(t) = \intt \theta_s^2 ds$, on both sides of the above equation we have 
\beo
d \Bl X_t^2 e^{F_1(t)} \Br &=& e^{F_1(t)}\theta_t^2 dt + e^{F_1(t)}X_t \theta_t dW_t \\
\Rightarrow X_t^{2} &=& e^{-F_1(t)} [X_0^2  + \intt e^{F_1(s)} \theta_s^2 ds + 2 \intt e^{F_1(s)} X_s \theta_s dW_s]\\
&=& e^{-F_1(t)}[ X_0^2 + \intt d(e^{F_1(s)}) + 2 \intt e^{F_1(s)} X_s \theta_s dW_s]\\
&=& e^{-F_1(t)}[X_0^2 + e^{F_1(t)}-1 + 2 \intt e^{F_1(s)} X_s\theta_s dW_s] \\
&=& X_0^2 e^{-F_1(t)} + 1 - e^{-F_1(t)} + 2 \intt e^{F_1(s)- F_1(t)} X_s \theta_s dW_s\\
\Rightarrow E(X_t^{2})&=& E(e^{-F_1(t)}) E(X_0^2) + 1 - E(e^{-F_1(t)}) \n \\
&+& 2 E \Bl e^{-F_1(t)} \intt e^{F_1(s)} X_s \theta_s dW_s \Br .
\eeo
From the proof of Lemma \ref{martingale lemma 2} we have that the third expectation is zero (by substituting $m = 1$).
Therefore 
\be
E(X_t^2) &=&  E(e^{-F_1(t)}) E(X_0^2) +1 -  E(e^{-F_1(t)})  \n \\
\Rightarrow \lim_{t \rightarrow \infty} E(X_t^2) &=& E(X_0^2) \lim_{t \rightarrow \infty} E(e^{-F_1(t)}) \label{EXt2} \\
&+& 1 -  \lim_{t \rightarrow \infty}E(e^{-F_1(t)})  \label{Xt2}
\ee
Now $\limtin E(e^{-F_k(t)}) = 0$ by Lemma \ref{Fkt}. Therefore, $$ \lim_{t \rightarrow \infty} E(X_t^2)= 1  \ \ \mbox{from (\ref{Xt2})}.$$

\item 
Assume this holds for $1 \le m \le k-1$, i.e.,
$$\lim_{t \rightarrow \infty} E(X_t^{2m}) =\frac{(2m)!}{2^m m!} \ \  \mbox{for} \ 1 \le m \le k-1.$$
\item 
From It\^{o}'s lemma applied to $X_t^{2k}$
\beo
dX_t^{2k} &=& \Bl-kX_t^{2k}\theta_t^2  + k(2k-1) X_t^{2k-2} \theta_t^2 \Br dt + 2k X_t^{2k-1} \theta_t dW_t .
\eeo
Multiplying with the integrating factor \ $e^{F_k(t)}$ on both sides of the above equation and rearranging we have that 
\be
d\Bl X_t^{2k} e^{F_k(t)} \Br &=& k(2k-1) e^{F_k(t)} X_t^{2k-2} \theta_t^2 dt + 2k e^{F_k(t)} X_t^{2k-1} \theta_t dW_t \n \\
\Rightarrow X_t^{2k} &=& e^{-F_k(t)}[ X_0^{2k} + (2k-1) \intt k e^{F_k(s)} X_s^{2k-2} \theta_s^2 ds \n \\
&+& 2k \intt e^{F_k(s)} X_s^{2k-1} \theta_s dW_s\n \\
\Rightarrow E(X_t^{2k}) &=& E(e^{-F_k(t)}) E(X_0^{2k}) \n \\
&+& (2k-1) E(\intt k e^{-F_k(t)} e^{F_k(s)} X_s^{2k-2} \theta_s^2 ds) \n \\
&+ & E \Bl 2 e^{-F_k(t)} \intt  ke^{F_k(s)} X_s^{2k-1} \theta_s  dW_s \Br .\label{XtAt}
\ee
We have proved in Lemma \ref{martingale lemma 2} that the third expectation in the RHS of (\ref{XtAt}) is zero (by substituting m = k). Writing  
\be 
A_{k,2m-2}(t) &:=& E( e^{-F_k(t)} k \intt e^{F_k(s)} X_s^{2m-2} \theta_s^2 ds)  \n \\
&=& E (e^{-F_k(t)} \intt X_s^{2m-2} d(e^{F_k(s)})), \ \ \mbox{for} \ 1 \le m \le k \label{Ak2m-2}
\ee
we have,
\be
 E(X_t^{2k}) &=&     E(X_0^{2k}) E(e^{-F_k(t)}) \n \\
&+& (2k-1)   A_{k,2k-2}(t)  \label{Xt2k}.
\ee
\vspace{-15pt}
Now by the integration by parts we have, 
\beo
\intt X_s^{2m}d(e^{F_k(s)}) &=& X_t^{2m}e^{F_k(t)} - X_0^{2m} - \intt e^{F_k(s)} d(X_s^{2m}) \n \\
&=& X_t^{2m} e^{F_k(t)} - X_0^{2m} 
- \intt e^{F_k(s)} \Bl (-m X_s^{2m} \theta_s^2 \n \\
&+& m(2m-1) X_s^{2m-2} \theta_s^2) ds 
+ \intt 2m X_s^{2m-1} \theta_s dW_s \Br ,
\eeo
using $$ d X_t^{2m} = - m X_t^{2m} \theta_t^2 dt + m(2m-1) X_t^{2m-2} \theta_t^2 dt  + 2m X_t^{2m-1} \theta_t d W_t.$$ 
Therefore multiplying by $e^{-F_k(t)}$ on both sides of the above equation we have 
\beo
e^{-F_k(t)} \intt X_s^{2m} d(e^{F_k(s)}) &=& 
e^{-F_k(t)} \intt k \theta_s^2 e^{F_k(s)}X_s^{2m} ds \n \\
&=&  X_t^{2m} - X_0^{2m} e^{-F_k(t)} \\
&+& e^{-F_k(t)} \intt m e^{F_k(s)} X_s^{2m} \theta_s^2 ds \n \\
&-& e^{-F_k(t)} \intt m (2m-1) e^{F_k(s)}X_s^{2m-2} \theta_s^2 ds \n \\ 
&+& 2m e^{-F_k(t)} \intt  X_s^{2m-1} e^{F_k(s)} \theta_s d W_s . \n
\eeo
Taking expectations on both sides and recalling the definition of $A_{k,2m}(t)$ from (\ref{Ak2m-2}) we have ,
\be
A_{k,2m}(t) &=& E(X_t^{2m}) - E(e^{-F_k(t)}) E(X_0^{2m}) + \frac{m}{k}A_{k,2m}(t) \n \\
&-& \frac{m(2m-1)}{k} A_{k,2m-2}(t) +  0. 
\ee
That the last expectation is zero follows from Lemma \ref{martingale lemma 2}. This implies that 
\be
(1 - \frac{m}{k}) A_{k,2m}(t) &=& E(X_t^{2m}) - E(e^{-F_k(t)})E(X_0^{2m}) \n \\
&-& \frac{m(2m-1)}{k} A_{k,2m-2}(t)
\label{Ak2m}.
\ee
Now,
\be
A_{k,0}(t) &=& E(e^{-F_k(t)} \intt k \theta_s^2 e^{F_k(s)} ds) \n \\
&=& E(e^{-F_k(t)} \intt d(e^{F_k(s)})) = 1 - e^{-F_k(t)}. \label{Ak2}
\ee
Define $B_{k,2m}= \lim_{t \rightarrow \infty} A_{k,2m}(t)$ (when the limit exists). Taking limits as $t \to \infty$ on both sides of (\ref{Ak2}) and applying Lemma \ref{Fkt} we get:
\be
B_{k,0} = 1 - \lim_{t \rightarrow \infty} e^{-F_k(t)} = 1. \label{Bko}
\ee
Hence $\lim \limits_ {t \to \infty} A_{k,2m}(t)$ exists for $m=0$.\\
Taking $m = 1,2,3,\ldots,k-1$ in (\ref{Ak2m}) we get that $\lim \limits_ {t \to \infty}A_{k,2m}(t)$ exists, since 
\be
(1 - \frac{m}{k}) \lim _{t \to \infty} A_{k,2m}(t) &=& \lim_{t \to \infty} E(X_t^{2m}) - \frac{m(2m-1)}{k} \lim_{t \to \infty} A_{k,2m-2}(t)\n \\
\R B_{k,2m} &=& \frac{k}{k-m} \lim_{t \to \infty} E(X_t^{2m}) - \frac{m(2m-1)}{k-m} B_{k,2m-2}. \n \\
\label{B1k}
\ee

Substituting different values of $m=0,1,2,\ldots,k-1$ in (\ref{B1k}) and applying induction hypothesis, that $\limtin E(X_t^{2m}) = \frac{(2m)!}{2^m m!}$, for $0\le m \le k-1$, we get:
\beo
B_{k,0}&=& 1\\
B_{k,2} &=& \frac{k}{k-1}1-\frac{1}{k-1}1 = 1  \\
B_{k,4}&=& \frac{k}{k-2}3 - \frac{2.3}{k-2}1 = 3 \n \\
B_{k,6}&=& \frac{k}{k-3}5.3-\frac{3.5}{k-3}3=5.3 \n \\ 
B_{k,8}&=& \frac{k}{k-4}7.5.3 - \frac{4.7}{k-4}5.3 = 7.5.3 \n \\
\ldots\\
B_{k,2k-2}&=& k (2k-3)(2k-5)\ldots 3.1 - (k-1)(2k-3)B_{k,2k-4}\\
&=& k(2k-3)(2k-5)\ldots3.1 -(k-1)(2k-3) \ (2k-5)\ldots3.1  \n \\
&=&(2k-3)(2k-5)\ldots3.1(k-k+1)\\
&=&\frac{(2k-2)!}{2^{k-1}(k-1)!} .
\eeo
\vspace{-10pt}
Therefore applying Lemma \ref{Fkt} to Equation (\ref{Xt2k}) :
\be
\lim_{t \rightarrow \infty} E(X_t^{2k}) &=& (2k-1) B_{k,2k-2}\n \\
&=& (2k-1) \frac{(2k-2)!}{2^{k-1}(k-1)!} = \frac{2k(2k-1)!}{2^k k!} \n \\
&=& \frac{(2k)!}{2^k k!}. \label{Xteven} 
\ee
\end{enumerate}

\textbf{Odd moments}
\begin{enumerate} 

\item
To find the odd moments of $X_t$ we perform similar procedure as above. We have 
\be
d X_t = -X_t \frac{\theta_t^2}{2} dt + \theta_t dW_t  \label{Xt}
\ee
Define $G_k(t)= \frac {2k+1}{2} \intt \theta_s^2 ds, k \in \mathbb{N} \cup \{0 \}$. Multiply by the integrating factor \ $e^{G_0(t)}$ \ 
on both sides of (\ref{Xt}) and rearrange to get 
\be
d(e^{G_0(t)} X_t) &=& e^{G_0(t)} \theta_t d W_t \n \\
\Rightarrow X_t &=& X_0 e^{-G_0(t)} + e^{-G_0(t)} \intt e^{G_0(s)} \theta_s dW_s \n \\
\Rightarrow E(X_t) &=& E(X_0) E(e^{-G_0(t)}) + E \Bigl( e^{-G_0(t)} \intt e^{G_0(s)} \theta_s d W_s \Bigr). \label{EXt}
\ee
From Lemma \ref{Fkt} we have 
$$\lim_{t \rightarrow \infty} E(e^{-G_0(t)}) = 0.
$$

Therefore from (\ref{EXt}) we have 
$$  \lim_{t \rightarrow \infty} E(X_t) = 0.$$
\item Let $k \ge 1$ be any positive integer. Assume that
$$\lim_{t \rightarrow \infty} E(X_t^{2m-1})=0 \ \mbox{where}\ \  m=1,2,\ldots, k.  $$
\item Applying It\^{o}'s lemma to $X_t^{2k+1}$ we get 
\be
d X_t^{2k+1}&=& (2k+1) X_t^{2k} d X_t + \frac{1}{2}(2k+1)2k X_t^{2k-1} \theta_t^2 dt \n \\
&=& (2k+1) X_t^{2k} \Bigl( -X_t \frac{\theta_t^2}{2} dt + \theta_t d W_t \Bigr) + (2k+1)k \theta_t^2 X_t^{2k-1} dt \n \\
&=& \Bigl( - \frac{1}{2}(2k+1) X_t^{2k+1} \theta_t^2 \n \\
&+& (2k+1) k X_t^{2k-1} \theta_t  \Bigr)dt + (2k+1) \theta_tX_t^{2k}dW_t. \label{X_2kp1}
\ee
Multiplying by the integrating factor $e^{G_k(t)}$ on both sides of (\ref{X_2kp1}) and rearranging we get:
\be 
d \Bigl(  X_t^{2k+1} e^{G_k(t)} \Bigr) &=& k(2k+1) e^{G_k(t)} \theta_t^2 X_t^{2k-1} dt + (2k+1) e^{G_k(t)} \theta_t X_t^{2k} dW_t \n \\
\Rightarrow X_t^{2k+1} &=& e^{-G_k(t)} \Bigl[  X_0^{2k+1} +k(2k+1) \intt e^{G_k(s)} \theta_s^2 X_s^{2k-1}ds \n \\
&+& (2k+1) \intt e^{G_k(s)} \theta_s X_s^{2k} dW_s \Bigr]. \n
\ee
Thus
\be 
E(X_t^{2k+1}) &=& E(e^{-G_k(t)}) E(X_0^{2k+1}) \n \\
&+& E\Bigl( k(2k+1) e^{-G_k(t)}\intt e^{G_k(s)}  \theta_s^2 X_s^{2k-1} ds \Bigr) \n \\
&+& (2k+1) E \Bl e^{-G_k(t)} \intt  e^{G_k(s)} X_s^{2k} \theta_s dW_s \Br. \label{Xt2kp1}
\ee
From Lemma \ref{rmk mart} we have the third expectation is zero. That is 
$$ E \Bl e^{-G_k(t)} \intt e^{G_k(s)} X_s^{2k} \theta_s d W_s \Br = 0.$$ 

Defining 
\be 
C_{k,2m-1}(t) &:=& E\Bigl( k(2k+1) e^{-G_k(t)} \intt e^{G_k(s)}\theta_s^2 X_s^{2m-1} ds  \Bigr)\n \\
&=& E\Bigl( 2k e^{-G_k(t)} \intt X_s^{2m-1}d(e^{G_k(s)}) \Bigr). \label{Ck2m-1}
\ee
We have from (\ref{Xt2kp1}). 
\be
 E(X_t^{2k+1}) &=&  E(e^{-G_k(t)}) E(X_0^{2k+1})+  C_{k,2k-1}(t).  \label{Xtodd}
\ee
Now by integration by parts 
\be	
\intt X_s^{2m-1} d(e^{G_k(s)})=  X_t^{2m-1} e^{G_k(t)} - X_0^{2m-1} - \intt e^{G_k(s)} d(X_s^{2m-1}). \n 
\ee
Applying It\^{o}'s lemma to $X_t^{2m-1}$ we have 
\be
dX_t^{2m-1} &=&  (2m-1) X_t^{2m-2} dX_t + (2m-1)(m-1)X_t^{2m-3} \theta_t^2 dt  \n \\
&=& - \frac{2m-1}{2} X_t^{2m-1} \theta_t^2 dt + (m-1)(2m-1) X_t^{2m-3} \theta_t^2 dt + (2m-1) X_t^{2m-2} \theta_t d W_t. \n 
\ee
Substituting in the above equation we have 
\be 
\intt X_s^{2m-1} d (e^{G_k(s)}) &=& X_t^{2m-1} e^{G_k(t)} - X_0^{2m-1} + \intt (2m-1) e^{G_k(s)} X_s^{2m-1}\frac{\theta_s^2}{2} ds \n \\
&-& (2m-1)(m-1) \intt e^{G_k(s)} X_s^{2m-3} \theta_s^2 ds - \intt (2m-1) e^{G_k(s)} X_s^{2m-2} \theta_s dW_s.  \n \\
\label{inter1}
\ee
Multiplying both sides by $e^{-G_k(t)}$, taking expectations in  (\ref{inter1}) and recalling the definition of $C_{k,2m-1}$  from (\ref{Ck2m-1}) we have 
\be
C_{k,2m-1}(t) &=& E\Bl 2k e^{-G_k(t)} \intt X_s^{2m-1} d(e^{G_k(s)})\Br \n \\
&=& 2k E(X_t^{2m-1}) - 2k E(e^{-G_k(t)}) E(X_0^{2m-1}) + \frac{2m-1}{(2k+1)} C_{k,2m-1}(t) \n \\
&-& \frac{(2m-1)(2m-2)}{(2k+1)} C_{k,2m-3}(t) - (2m-1) E  \Bl e^{-G_k(t)} \intt e^{G_k(s)} X_s^{2m-2} \theta_s dW_s \Br. \n \\
 \label{ckm}
\ee
Now by Lemma \ref{rmk mart} where it is shown that $$E \Bl e^{-G_k(t)} \intt \theta_s e^{G_k(s)} X_s^{2m} dW_s \Br = 0, \ \mbox{for} \ 0 \le m \le k-1,$$ we have that the third expectation is zero. 
Now  
\beo
C_{k,1}(t)&=&  2k E \Bl e^{-G_k(t)}\intt X_s d(e^{G_k(s)})\Br \n \\
&=& 2k E \Bl e^{-G_k(t)} \Bl X_t e^{G_k(t)} - X_0 e^{G_k(0)} - \intt e^{G_k(s)} d X_s \Br \Br \n \\
&=& 2k E \Bl X_t - X_0 e^{-G_k(t)} - e^{-G_k(t)}\intt e^{G_k(s)} dX_s \Br \ \mbox{since $G_k(0)=0$.}
\eeo
From the SDE of $X_t$ we have 
\beo
C_{k,1}(t) &=& 2k E \Bl X_t - X_0e^{-G_k(t)} - e^{-G_k(t)} \intt e^{G_k(s)} \Bl -X_s\frac{\theta_s^2}{2} ds + \theta_s dW_s \Br \Br \n \\
&=& 2k \Bl E(X_t) -E(X_0 e^{-G_k(t)}) + \frac{1}{2} E \Bl e^{-G_k(t)}\intt e^{G_k(s)} \theta_s^2 X_s ds \Br \n \\
&-& E \Bl e^{-G_k(t)} \intt e^{G_k(s)} \theta_s d W_s \Br \Br \n \\
&=& 2k \Bl E(X_t) - E(X_0e^{-G_k(t)}) \Br + \frac{1}{2k+1} C_{k,1}(t) + 0,\n 
\eeo
since $E \Bl e^{-G_k(t)} \intt e^{G_k(s)} \theta_s dW_s\Br =0$ from Lemma \ref{rmk mart} (by substituting m = 0). Therefore
\be
(1 - \frac{1}{2k+1})C_{k,1}(t) &=& 2k \Bl E(X_t) -  E(X_0 e^{-G_k(t)})\Br. \label{Ck1}
\ee
Now we have proved that $\limtin E(X_t) = 0 = \limtin E(e^{-G_k(t)})$. Defining $D_{k,m} = \limtin C_{k,m}(t)$ for $m = 1,3,5,\ldots,2k-1$, wherever it exists, we have from (\ref{Ck1})
$$
D_{k,1} = 0.
$$
From (\ref{ckm}) we have 
\be
(1  - \frac{2m-1}{2k+1})C_{k,2m-1}(t)&=& 2k E(X_t^{2m-1}) - 2k E(e^{-G_k(t)}) E(X_0^{2m-1}) \n \\
&-& \frac{(2m-1)(2m-2)}{(2k+1)} C_{k,2m-3}(t).\label{Ck2m}
\ee
By induction hypothesis $\limtin E(X_t^{2m-1}) = 0$ for $m = 1,2,\ldots,k$.  
Since $D_{k,1} = 0$ from (\ref{Ck2m}) we have by iteration $\limtin C_{k,2m-1}(t)$ exists and equals to 0 for $m =1,2,\ldots,k$, i.e. 
$$ D_{k,j} = 0 \ \ \mbox{for $j = 1,3,\ldots,2k-1$}. $$

Therefore, from \ref{Xt2kp1} we have that 
\be
\limtin E(X_{t}^{2k+1}) = 0. \label{xtodd}
\ee
\end{enumerate}

Thus combining (\ref{Xteven}) and (\ref{Xtodd}) we see that the limiting moments of $\{ X_s\}$
matches with that of a $N(0,1)$ distribution. Since the limiting distribution admits a smooth density, invoking uniqueness of moment generating function we can infer that the limiting distribution of $\{ X_s \}$ is $N(0,1)$. This completes the proof of Theorem \ref{main Theorem}.
$\hfill{\blacksquare}$

\begin{rmk}
\label{lower bound of theta}
From (\ref{eta form}) we have $\theta_t$ satisfying the equation
\be
\theta_t &=& \frac{e^{qt}}{\eta_0 + \frac{1}{\sqrt{2 \pi}} \intt e^{qs}|X_s| ds } \n \\
\Rightarrow \theta_t^2 &=& \frac{e^{2qt}}{\Bl\eta_0 + \frac{1}{\sqrt{2 \pi}} \intt e^{qs}|X_s| ds  \Br^2} \n \\
&\ge& \frac{e^{2qt}}{2 \Bl \eta_0^2 + \frac{1}{2 \pi} (\intt e^{qs}|X_s|ds)^2 \Br}\n \\
&=& \frac{e^{2qt}}{2 \eta_0^2 + \frac{(e^{qt} -1)^2}{ \pi q^2} (\intt \frac{q}{e^{qt}-1} e^{qs} |X_s| ds)^2}\n \\
&\ge& \frac{e^{2qt}}{2 \eta_0^2 + \frac{e^{qt}-1}{\pi q} \intt e^{qs} |X_s|^2 ds}, \label{thetasq2}
\ee
where the last inequality follows from the fact that 
\beo
(\intt \frac{q}{e^{qt}-1} e^{qs}|X_s| ds)^2 &\le& \intt \frac {q}{e^{qt}-1} e^{qs} |X_s|^2 ds.
\eeo
This is true by the Jensen's  inequality 
\beo
(E|X_s|)^2 &\le& E(|X_s|^2),
\eeo
with the expectation computed with respect to the density $f(x) = \frac{q}{e^{qt}-1} e^{qx}, 0 < x < t,$ for any $t>0$.
Therefore
\beo
E(\theta_t^2)&\ge& \frac{e^{2qt}}{2 E(\eta_0^2) + \frac{e^{qt}-1}{\pi q} \intt e^{qs} E(X_s^2) ds} \\
&\ge& \frac{e^{2qt}}{2 E(\eta_0^2) + \frac{(e^{qt}-1)^2(1 + E(X_0^2))}{\pi q^2}},
\eeo
where the last inequality follows from (\ref{EXt2}) that 
$$ E(X_t^2) \le 1 + E(X_0^2) \ \forall t > 0.$$
Therefore
\beo
\liminf_{t \to \infty} E(\theta_t^2) &\ge & \liminf_{t \to \infty} \frac{e^{2qt}}{2 E(\eta_0^2) + \frac{(e^{qt} - 1)^2 (1+E(X_0^2))}{\pi q^2}} = \frac{\pi q^2}{1+E(X_0^2)}  \\
\Rightarrow  \liminf_{t \to \infty} \frac{1}{t} \intt E(\theta_u^2) du &\ge&  \frac{\pi q^2}{1+E(X_0^2)}, \ \ \mbox{by Fatou's lemma}.
\eeo
In particular if $X_0=0$ almost surely, then 
\beo
\liminf_{t \to \infty} \frac{1}{t} \intt E(\theta_u^2) du &\ge& \pi q^2.
\eeo
This gives a lower bound to the growth of $\theta_t$.
\end{rmk}

\begin{rmk}
\label{rates of convergence}
\textbf{Rates of convergence of Adaptive and Standard MCMC:} Recalling the SDE for AMCMC for Normal target density for $X_t$ is given as:
\beo
d X_t &=& - X_t \frac{\theta_t^2}{2} dt + \theta_t d W_t.
\eeo
Multiplying by the integrating factor and performing the usual operations we get:
\be
E(X_t)&=& E(X_0 e^{-G_t}), \ \mbox{where $G_t = \intt \frac{\theta_s^2}{2}ds $.} \label{xt1}
\ee
Similar equation for the SMCMC $Y_t$ is:
\beo
d Y_t &=& -Y_t \frac{\theta_0^2}{2} dt + \theta_0 dW_t.
\eeo
Applying similar computations we get 
\beo
E(Y_t)&=& E(Y_0 e^{-\tilde{G}_t}), \ \mbox{where $\tilde{G}_t = \intt \frac{\theta_0^2}{2}ds = {\frac{\theta_0^2}{2}t}$.}
\eeo
Similar computation with $X_t^2$ will give (see the proof of Lemma \ref{even odd})
\be
E(X_t^2) &=& E(X_0^2 e^{-G(t)}) + 1 - E(e^{-G(t)}). \label{xt2}
\ee
It is therefore clear from the Equations (\ref{xt1}) and 
(\ref{xt2}) that the quantity regulating the speed to 
convergence is $G(t)$ (or $\tilde{G}(t)$). The faster $G(t)$ 
(or $\tilde{G}(t)$) goes to $\infty$, the faster the process 
converges to its invariant distribution (which is standard 
Normal in this case). For the diffusion defined by the SDE corresponding to the SMCMC the 
rate of convergence to its stationary distribution is 
exponentially fast in $t$. For the AMCMC it depends on the 
behaviour of $\intt \theta_s^2 ds$. We have shown in Lemma 
\ref{timeave lem} that $\limsup \limits \limits_{t \to 
\infty} \frac{1}{t} \intt E(\theta_u^2) du < C < \infty$ for 
any $k \in \mathbb{N}$ when the target distribution is 
standard Normal. Combining this with Remark \ref{lower bound 
of theta} we find that the rate of convergence of the process defined by the SDE for 
the AMCMC to its stationary distribution is exponentially 
fast with exponent is linear in $t$. Thus the comparison 
between the rate of convergence of the processes defined by the SDEs for the  AMCMC and 
SMCMC to their corresponding stationary distribution will 
depend on the lower bound $\pi q^2$ and the upper bound $C$ 
(as in Lemma \ref{timeave lem} for $k=4$) and $\theta_0^2$. 
\tr{If the bound can be obtained in the almost sure sense, 
and not in the $L_1$ sense then it might be possible to 
directly compare SMCMC and AMCMC.}
\end{rmk}
\begin{rmk}
It is true that for the discrete time SMCMC, higher value of $\theta_0$ will delay convergence to stationarity of the chain. 
However, it is somewhat misleading that the diffusion process corresponding to the SMCMC converges faster to its
stationary distribution for higher value of $\theta_0$.  
For the AMCMC situation is quite different. The simulations in Figures 3.1 and 3.2 show that the trajectories of $\theta_t$ converge for large values of $t$. This is in tune to our theoretical findings that for a standard Normal target with standard Normal proposals, the time average moments of $\theta_t$ are bounded. Since this happens for any starting value of $\theta_0$, we recommend that this limiting value (or variable) should be used for selecting the optimal value of $\theta_0$. One should run the AMCMC sufficiently long, till the point where $\theta_t$ changes no further or varies in a narrow range. From that point onwards one should keep the level of $\theta_t$ same (any point in the narrow interval) and run a simple SMCMC. 
\end{rmk}

\section{Conclusion}
\label{conclusion}
Verifying  the conditions of Roberts \textit{et al}. 
 for checking the ergodicity of an AMCMC can 
sometimes prove to be difficult. In Basak and Biswas \citep{Basak}, 
we considered an AMCMC with the proposal
 kernel dependent on the previously generated sample and an arbitrary target distribution. There we performed
 a diffusion approximation technique to look at the continuous time version of the discrete chain. In this
 paper we narrowed down to the case where the target distribution is standard Normal. We investigate whether 
the invariant distribution of the diffusion is indeed the target distribution. It turns out that the resulting 
diffusion (which although singular) admits a unique invariant distribution. Then computing the limiting moments (both even and odd)
of $X_t$ we identify the limiting distribution to be $N(0,1)$.\\
The techniques applied here are specific only when the target distribution is Normal. Different methodologies may be needed to extend these results to other target distributions, where an identification of the limiting moments may not be possible.
 Also more choices of the proposal distribution can be made, where the kernel is dependent on a finite 
(or possibly infinite) past. We plan to take up these issues in our future work.

{
}

\end{document}